\documentclass{amsart}
\usepackage{amssymb}
\usepackage{amscd}
\DeclareFontFamily{OMS}{cmsy}{%
\fontdimen16\font=3pt
\fontdimen17\font=3pt}
\makeatletter
\renewcommand{\subsection}{\@startsection{subsection}{2}{\z@}%
{\baselineskip}{0.5\baselineskip}{\bfseries}}
\makeatother
\def\dj{d\kern-.30em\raise1.25ex\vbox{\hrule width .3em height .03em}}
\def\Dj{D\rlap{\kern-.70em\raise0.75ex
\vbox{\hrule width .3em height .03em}}}
\newtheorem{thm}{Theorem}
\newtheorem{pro}[thm]{Proposition}
\newtheorem{lem}[thm]{Lemma}
\theoremstyle{definition}
\newtheorem{defn}{Definition}
\numberwithin{equation}{section}
\newenvironment{pf}{\proof[\proofname]}{\endproof}
\def\adjr{\varkappa}
\def\dM{\mbox{$\smash{\vphantom{d}^{\raise0.4ex\hbox{$\scriptscriptstyle M$}}%
\!\!\!d}$}}
\def\V{\Bbb{V}}
\def\v{\varkappa}
\def\bla#1{$(${\it #1\/{}}$)$}
\def\ll{L}
\def\SO{\mathrm{SO}}
\def\U{\mathrm{U}}
\def\M{\mathrm{M}}
\def\D{\nabla}
\def\VB{\cal{H}}
\def\restr{\restriction}
\def\gen{\mathrm{gen}}
\def\cal{\mathcal}
\def\Bbb{\mathbb}
\def\e{\epsilon}
\def\k{\kappa}
\def\ad{\mathrm{ad}}
\def\id{\mathrm{id}}
\def\grten{\mathbin{\widehat{\otimes}}}
\def\adj{\varpi}
\def\pre{\mathfrak{der}}
\def\prv{\overrightarrow{\mathfrak{der}}}
\def\hor{\mathfrak{hor}}
\def\rig{\wp}
\def\inv{{i\!\hspace{0.8pt}n\!\hspace{0.6pt}v}}
\def\im{\mathrm{im}}
\def\Sum{{\displaystyle\sum}}
\def\L{\cal{L}}
\def\pHG{\varsigma}
\def\AH{\cal{A}'}
\def\phiH{\phi'}
\def\phiHG{\phi_*}
\def\piH{\pi'}
\def\kH{\k'}
\def\adjH{\adj'}
\def\eH{\e'}
\def\AHi{\cal{K}}
\def\adHG{\ad_*}
\begin{document}
\title[quantum principal bundles]{GENERAL FRAME STRUCTURES ON\\
QUANTUM PRINCIPAL BUNDLES}
\author{Mi\'co {\Dj}ur{\Dj}evi\'c}
\address{Instituto de Matematicas, 
UNAM, Area de la Investigacion Cientifica, 
Circuito Exterior, Ciudad Universitaria, M\'exico DF, cp 04510,
MEXICO}
\begin{abstract} A noncommutative-geometric generalization of the
classical formalism of frame bundles is developed, incorporating
into the theory of quantum principal bundles the concept
of the Levi-Civita connection. The construction of 
a natural differential calculus on quantum principal frame bundles is
presented, including the construction of
the associated differential calculus on the structure group. 
General torsion operators are defined and analyzed.
Illustrative examples are presented.
\end{abstract}
\maketitle

\section{Introduction}
\renewcommand{\thepage}{}

The formalism of principal bundles plays a central role in
the foundation of classical differential geometry.
However as far as basic differential-geometric structures are
concerned, the appropriate
conceptual framework is given by a more restrictive class of frame bundles.
Generally speaking, frame bundles are understandable as
(covering bundles) of the appropriate
subbundles of the principal bundle of linear frames
associated to the base manifold. 
These bundles intrinsically express the geometrical structure existing on the
base.

Various fundamental concepts and ideas of theoretical
physics are also most naturally and effectively
formulated and studied at the language of principal bundles. As a paradigmic
example, let us mention classical and quantum gauge theories.
The basic example of applying the formalism of frame bundles in
theoretical physics is given by general relativity theory. 

In this paper we are going to incorporate the classical theory of
principal frame bundles into the framework of non-commutative
differential geometry \cite{c}. All considerations are logically based
on a general theory of quantum principal bundles \cite{d1,d2} in which
the bundle and the base are quantum objects, and quantum groups
play the role of the structure groups. We shall consider
here only bundles with compact structure groups. Our constructions
are motivated by properties of classical frame structures with
orthonormal frames, which preassumes the existence of the metrics on
the base manifold. However, the developed formalism can be applied in the
general context, dealing with bundles with arbitrary structure groups. 

The paper is organized in the following way.
The next section is devoted to the definition of frame structures on
quantum principal bundles. In classical geometry, every framed bundle
involving orthonormal frames canonically determines the unique metric
connection with vanishing torsion. This is by definition the Levi-Civita
connection. It is also of the special interest to consider subbundles of
the orthonormal frame bundle, invariant under the Levi-Civita
connection.

The starting idea of the next section is to generalize the classical
concept of coordinate first-order forms \cite{kn} which determine the
association of orthonormal frames of the tangent spaces of the base
manifold, to the points of the corresponding fibers of the bundle.
This line of thinking was also suggested in \cite{ph}, in a conceptually
different context. 
Formally, we shall start from a vector space $\V$ representing abstract
quantum first-order coordinate horizontal forms, which is equipped with
the appropriate additional structure, equivalent to specifying a
bicovariant *-bimodule $\Psi$ over \cite{w2} the structure quantum
group $G$. The space $\V$ is then re-interpretable as the left-invariant
part of $\Psi$. In particular, there exists the intrinsic braid operator
$\tau\colon \V\otimes \V\rightarrow \V\otimes \V$. Using this operator, we
can construct the counterpart of the corresponding exterior algebra
$\V^\wedge$, requiring the $\tau$-antimultiplicativity between the
elements of $\V$. In our context this corresponds to the subalgebra
of horizontal forms generated by coordinate forms. Finally,
the complete algebra $\hor_P$ representing horizontal forms can be
constructed by taking the natural cross product between $\V^\wedge$
and the algebra $\cal{B}$ representing the (appropriate functions on the)
quantum principal bundle $P$. The concept of a frame structure is
completed by
introducing the analogs of Levi-Civita connections. In the formalism,
they will be represented by certain antiderivations
$\D\colon\hor_P\rightarrow\hor_P$, playing the role of the covariant
derivative maps of actual Levi-Civita connections. 

\renewcommand{\thepage}{\arabic{page}}
Starting from a framed quantum principal bundle $P$, and applying a
general construction presented in \cite{d3}, it is possible to
construct in a natural manner the whole differential calculus
$\Omega(P)$ on $P$, including a natural bicovariant first-order
*-calculus $\Gamma$ on $G$. This problematics is discussed in
Section~3. In particular, it turns out that $\hor_P$ coincides with the
algebra of horizontal forms associated to $\Omega(P)$, as defined
in the framework of the general theory \cite{d2}. Moreover, regular
(and multiplicative) connections \cite{d2} on $P$ are in a natural
correspondence with special antiderivations $D\colon\hor_P\rightarrow
\hor_P$. These maps play the role of the corresponding
covariant derivative maps. In particular, the construction justifies the
mentioned interpretation of the frame structure. 

In Section~4 we shall introduce and analyze counterparts of torsion
operators, associated to arbitrary connections on $P$. Section~5 gives
the explicit `coordinate' description of covariant derivative maps
$D$, in terms of the analogs of associated coordinate horizontal vector
fields.

Finally, in Section~6 some interesting examples are collected, and some
concluding remarks are made. A large class of examples of quantum
frame bundles comes from the appropriate quantum homogeneous spaces \cite{d2}
and bundles associated to them. It turns out that if $G$ is a subgroup of
a quantum group $H$, and if $H$ is equipped with a special
first-order differential calculus, then every principal
$H$-bundle becomes a framed principal $G$-bundle, after
restricting the action map to the subgroup $G$. A theory of quantum frame
bundles with classical structure groups \cite{d-frm} will be also considered
from the point of view of the general formalism. A particular attention
will be given to the framed quantum line bundles, where it is possible
to write down simple expressions for all basic entities figuring
in the game.

\section{Framed Quantum Principal Bundles}

Let $G$ be a compact matrix quantum group \cite{w2},
represented by a Hopf *-algebra $\cal{A}$. The elements of
$\cal{A}$ are interpreted as  
`polynomial functions' on $G$. The coproduct, counit and the
antipode map will be
denoted by $\phi$, $\e$ and $\k$ respectively.

Let us consider a bicovariant \cite{w2} bimodule $\Psi$ over $G$, equipped
with the left/right action
maps $\ell_\Psi\colon\Psi\rightarrow\cal{A}\otimes\Psi$ and
$\rig_\Psi\colon\Psi\rightarrow\Psi\otimes\cal{A}$. There exists a
natural decomposition 
$\Psi\leftrightarrow\cal{A}\otimes \V$,
where $\V=\Psi_{\inv}$ is the corresponding left-invariant part.

The bicovariant bimodule structure on $\Psi$ is re-expressed via the right
action $\v =(\rig_\Psi{\restr}\V)\colon \V\rightarrow \V\otimes\cal{A}$
and the natural right $\cal{A}$-module structure
$\circ\colon \V\otimes\cal{A}\rightarrow \V$, where
$\vartheta\circ a=\k(a^{(1)})\vartheta a^{(2)}$.

Let us assume that $\Psi$ is also *-covariant. Then the space $\V$ is 
*-invariant and the following compatibility
conditions hold:
\begin{gather*}
\v *=(*\otimes *)\v \qquad (\vartheta\circ a)^*=\vartheta^*\circ\k(a)^*\\
\v (\vartheta\circ a)=\sum_k(\vartheta_k\circ a^{(2)})\otimes\k(a^{(1)})c_k
a^{(3)},
\end{gather*}
where $\Sum_k\vartheta_k\otimes c_k=\v (\vartheta)$.

Let us denote by $\tau\colon \V\otimes \V\rightarrow \V\otimes \V$
the canonical braid operator \cite{w2} associated to $\Psi$. It is
expressed via $\v $ and $\circ$ as follows
$$\tau(\eta\otimes\vartheta)=\sum_k\vartheta_k\otimes(\eta\circ c_k).$$
Let $\V^\wedge$ be the corresponding $\tau$-exterior algebra, obtained by
factorizing the tensor algebra
$\V^\otimes$ through the quadratic relations $\im(I+\tau)$. At this point
it is reasonable to assume that the braid operator $\tau$ is such that
$\ker(I+\tau)\neq\{0\}$. This condition 
ensures the nontriviality of the higher-order part of $\V^\wedge$.

Let us also assume that the space $\V$ is equipped with a scalar product, such
that the action $\v \colon\V\rightarrow\V\otimes\cal{A}$ is unitary. Let us
fix an orthonormal basis $\{\theta_1,\dots,\theta_n\}$ in the space $\V$.
In this basis, the action is described by a unitary matrix
$\{u_{ij}\}$, given by $\v (\theta_i)=\Sum_j\theta_j\otimes u_{ji}$.

The *-structure, $\v $ and $\circ$ are naturally extendible to
$\V^{\otimes,\wedge}$, by requiring
\begin{gather*}
(\eta\vartheta)\circ a=(\eta\circ a^{(1)})(\vartheta\circ a^{(2)})\qquad
1\circ a=\e(a)1\\
(\eta\vartheta)^*=(-1)^{\partial\eta\partial\vartheta}\vartheta^*\eta^*\qquad
\v (\eta\vartheta)=\v (\eta)\v (\vartheta).
\end{gather*}
The extended maps define
bicovariant graded *-algebras $\Psi^{\otimes,\wedge}\leftrightarrow
\cal{A}\otimes \V^{\otimes,\wedge}$.

Let $P=(\cal{B},i,F)$ be a quantum principal $G$-bundle \cite{d2}
over a quantum space $M$. By definition, $\cal{B}$ is a *-algebra
representing $P$ at the level of 
quantum spaces, while $F\colon\cal{B}\rightarrow\cal{B}\otimes\cal{A}$ is
a *-homomorphism playing the role of the dualized right action of $G$ on $P$.
Finally, $i\colon\cal{V}\rightarrow\cal{B}$ is the dualized projection
of $P$ on $M$, and $\cal{V}$ is the *-algebra representing $M$.
The image of $i$ coincides with the $F$-fixed-point
subalgebra of $\cal{B}$. Geometrically this means that
$M$ is interpretable as
the orbit space associated to $P$, relative to the action of the
structure group. We shall identify the elements of $\cal{V}$ with
their images in $\cal{B}$. 

\begin{lem} \bla{i} The formulas
\begin{align*}
(q\otimes\vartheta)(b\otimes\eta)&=\sum_kqb_k\otimes(\vartheta\circ c_k)
\eta\\
(b\otimes\vartheta)^*&=\sum_k b_k^*\otimes(\vartheta^*\circ c_k^*)
\end{align*}
where $\Sum_i b_i\otimes c_i=F(b)$,
determine a *-algebra structure on the graded vector space
$\hor_P=\cal{B}\otimes \V^\wedge$. In particular, $\cal{B}=\hor_P^0$ and
$\V^\wedge$ is a graded *-subalgebra of $\hor(P)$.

\smallskip\bla{ii}
There exists the unique *-homomorphism
$F^\wedge\colon\hor_P\rightarrow\hor_P\otimes\cal{A}$
extending actions $F$ and $\v \colon \V^\wedge
\rightarrow \V^\wedge\otimes\cal{A}$. The space $\hor_P$ is a graded
$\cal{A}$-comodule. \qed
\end{lem}

The above definition of the product in $\hor_P$, together with the
definition of $\V^\wedge$, implies commutation relations
$$ \vartheta\varphi=(-1)^{\partial\varphi\partial\vartheta}\varphi_k(
\vartheta\circ c_k),$$
where $\Sum_k \varphi_k\otimes c_k=F^\wedge(\varphi)$, while $\varphi\in
\hor_P$ and $\vartheta\in\V^\wedge$.

We pass to the definition of frame structures.

\begin{defn}\label{def:1} A frame structure on a quantum principal bundle $P$
relative to a bicovariant *-bimodule $\Psi$ is a first-order hermitian
antiderivation
$\D\colon\hor_P\rightarrow\hor_P$ such that

\smallskip
\bla{i} The following equalities hold
\begin{gather}
F^\wedge\D=(\D\otimes\id)F^\wedge\label{cov-cond}\\
\D(\VB)=\{0\},\label{lc-cond}
\end{gather}
where $\VB\subseteq\hor_P$ is the *-subalgebra generated by $\D(\cal{V})$ and
$\V$.

\smallskip
\bla{ii} There exist elements $b_{\alpha i}\in\cal{B}$ and $f_\alpha\in
\cal{V}$ such that 
\begin{align}
1\otimes\theta_i&=\sum_\alpha b_{\alpha i}\D(f_\alpha)\label{D-coor}\\
F(b_{\alpha i})&=\sum_jb_{\alpha j}\otimes u_{ji},
\end{align}
for each $i\in\{1,\dots, n\}$.
\end{defn}

Geometrically, the definition of a framed bundle incorporates the
idea that every point of the bundle gives rise to an orthonormal
system in the correponding tangent space, spanning a complement to
the vertical subspace. These complements form a special
connection on the bundle (actually the Levi-Civita connection), and the
map $\D$ is the corresponding covariant derivative.
An alternative possible geometrical interpretation
of the frame structures is that they represent
the appropriate `metrics' on the base space $M$.

\section{Differential Calculus On Framed Bundles}

In this section we describe the construction of the intrinsic
differential calculus, which can be associated to every quantum
principal bundle $P$ endowed with a frame structure. The construction
gives also a natural differential calculus on the structure group $G$.

Let $\Omega_M\subseteq\hor_P$ be the $F^\wedge$-fixed point subalgebra.
This is a graded *-subalgebra of $\hor_P$ and we have $\Omega_M^0=\cal{V}$.

Let us fix a frame structure $\D\colon\hor_P\rightarrow\hor_P$. 
The $F^\wedge$-covariance of
$\D$ implies that $\Omega_M$ is $\D$-invariant. In what follows, we shall
denote by
$\dM\colon\Omega_M\rightarrow\Omega_M$ the corresponding
restriction map.

\begin{lem} \bla{i} The algebra $\Omega_M$ is generated by elements of
the form $\{f,\dM(f)\}$, where $f\in\cal{V}$.

\smallskip\bla{ii} The map $\dM$ is a hermitian differential on
$\Omega_M$.
\end{lem}

\begin{pf} Inductively applying equation \eqref{D-coor} and the basic
commutation relation, it follows that the elements from
$\hor_P^m$ are of the form
$$\varphi=\sum b \prod_{i=1}^m \dM(f_i),$$
where $b\in\cal{B}$ and $f_i\in\cal{V}$. If $\varphi\in\Omega_M$ then it
is possible to assume that $b\in\cal{V}$, 
as follows from the fact that there exists a natural projection map
$p_0\colon\hor_P\rightarrow\Omega_M$ given by 
$$ p_0=(\id\otimes h)F^\wedge,$$
where $h\colon\cal{A}\rightarrow\Bbb{C}$ is the Haar measure \cite{w1}
on $G$. This proves \bla{i}.

To complete the proof, it is sufficient to apply \bla{i} and to
observe that the square of $\dM$ vanishes
on $\cal{V}$. It follows that $\dM^2=0$, globally.
\end{pf}

Let us consider a real affine space $\pre(P)$ consisting of hermitian
first-order antiderivations $D\colon\hor_P\rightarrow\hor_P$ intertwining
the map $F^\wedge$ and extending the differential $\dM$. By definition,
$\D\in\pre(P)$. 
Applying the results of \cite{d3}, it follows that for each $D\in\pre(P)$
there exist the unique linear map $\varrho_D\colon\cal{A}\rightarrow\hor_P$
such that
\begin{equation}
D^2(\varphi)=-\sum_k\varphi_k\varrho_D(c_k),
\end{equation}
where $\Sum_k\varphi_k\otimes c_k=F^\wedge(\varphi)$.

\begin{defn}
The map $\varrho_D$ is called the curvature of $D$.
\end{defn}

As explained in \cite{d3}, the curvature $\varrho_D$ also satisfies the
following identities:
\begin{align*}
F^\wedge\varrho_D&=(\varrho_D\otimes \id)\ad\\
D\varrho_D&=0\qquad \varrho_D(a)^*=-\varrho_D[\k(a)^*]\\
\varrho_D(a)\varphi&=\sum_k\varphi_k\varrho_D(ac_k)+\e(a)D^2(\varphi),
\end{align*}
for each $\varphi\in\hor_P$ and $a\in\cal{A}$. Here
$\ad\colon\cal{A}\rightarrow\cal{A}\otimes\cal{A}$ is the adjoint action
of $G$ on itself, given by $\ad(a)=a^{(2)}\otimes\k(a^{(1)})a^{(3)}$.

Furthermore, let us consider the real vector space $\prv(P)$ associated
to $\pre(P)$. The elements of $\prv(P)$ are interpretable as first-order
$F^\wedge$-covariant hermitian antiderivations on $\hor_P$ which
vanish on $\Omega_M$. Every such a map $E$ can be uniquely represented
in the form
\begin{equation}\label{E-k}
E(\varphi)=-(-1)^{\partial\varphi}\sum_k\varphi_k \chi(c_k)
\end{equation}
where $\chi=\chi_E\colon\cal{A}\rightarrow\hor_P$ is a first-order linear
map. According to \cite{d3}, the following equalities hold
\begin{gather}
F^\wedge\chi(a)=(\chi\otimes \id)\ad(a)
\label{E1}\\
\chi[\kappa(a)^*]=-\chi(a)^*\label{E2}\\
\chi(a)\varphi=\e(a)E(\varphi)+(-1)^{\partial\varphi}\sum_k\varphi_k\chi
(ac_k)\label{E3}
\end{gather}
for each $E\in\prv(P)$. Conversely, every first-order
linear map $\chi$ satisfying the above equalities determines an
antiderivation $E$, via equality~\eqref{E-k}.

\begin{lem}\label{lem:reduct} Every linear map $\chi=\chi_E$ satisfying
the above
mentioned properties is completely determined by its values on the matrix
elements of an arbitrary faithful representation of $G$. 
\end{lem}
\begin{pf}
Let $T\in\M_k(\Bbb{C})$ be a faithful unitary representation of $G$. This
means that the matrix elements of $T$ generate the *-algebra $\cal{A}$.
According to \eqref{E3} we have
\begin{equation}\label{kaT}
\chi(aT_{ij})=\sum_\alpha q_{\alpha i}\chi(a)p_{\alpha j},
\end{equation}
where $a\in\ker(\e)$ and the elements $q_{\alpha i},p_{\alpha i}\in\cal{B}$
satisfy
$$ F(p_{\alpha i})=\sum_j p_{\alpha j}\otimes T_{ji}\qquad
\sum_{\alpha} q_{\alpha i}p_{\alpha j}=\delta_{ij}1. $$
We can assume without a lack of generality that $T$ is
self-conjugate. Inductively applying \eqref{kaT} it follows that the values
of $\chi$ are algebraically expressible in terms of its values on matrix
elements $T_{ij}$.
\end{pf}

Equation \eqref{kaT} implies in particular the following generalization of
the corresponding classical antisymmetricity property
\begin{equation}\label{k-antisym}
\chi(T_{ij})+\sum_{k\alpha} q_{\alpha k}\chi(T_{ki}^*) p_{\alpha j}=0. 
\end{equation}

It is worth noticing that the same reasoning applies to the curvature maps
$\varrho_D$, and that the following connecting identity \cite{d3} holds
\begin{equation}\label{r-nu}
\varrho_{D+E}(a)=
\varrho_D
(a)+D\chi_E(a)+
\chi_E(a^{(1)})\chi_E(a^{(2)}). 
\end{equation}
The system of maps $\varrho_D$ and $\chi_E$ intrinsically
determines a bicovariant first-order *-calculus $\Gamma$ on $G$.
By definition,
this calculus is based on the right $\cal{A}$-ideal
$\cal{R}\subseteq\ker(\e)$
consisting of all elements anihilated by maps $\varrho_D$ and $\chi_E$.
As for all left-covariant differential structures, the left-invariant
part of $\Gamma$ is given by $\Gamma_{\inv}=\ker(\e)/\cal{R}$.

By definition, the maps $\varrho_D$ and $\chi_E$ are factorizable through
$\cal{R}$. In what follows we shall denote by the same symbols the
factorized maps $\varrho_D, \chi_E\colon\Gamma_{\inv}\rightarrow\hor_P$.
The following commutation relations hold
\begin{align}
\chi_E(\vartheta)\varphi&=(-1)^{\partial\varphi}
\sum_k\varphi_k\chi_E(\vartheta\circ c_k)\\
\varrho_D(\vartheta)\varphi&=\sum_k\varphi_k\varrho_D
(\vartheta\circ c_k)
\end{align}
for each $\vartheta\in \Gamma_{\inv}$
and $\varphi\in\hor_P$, where $\circ$ is the corresponding right 
$\cal{A}$-module structure on $\Gamma_{\inv}$. Furthermore, the maps
$\chi_E$ and $\varrho_D$ are hermitian, and intertwine $\adj$ and
$F^\wedge$. Here
$\adj\colon\Gamma_{\inv}\rightarrow\Gamma_{\inv}\otimes\cal{A}$ is the
induced adjoint action, explicitly given by
$$ \adj\pi=(\pi\otimes\id)\ad, $$
where $\pi\colon\cal{A}\rightarrow\Gamma_{\inv}$ is the canonical
projection map (playing the role of the germ map).

The described construction gives the minimal calculus $\Gamma$, compatible
with the whole space $\pre(P)$. The same construction can be performed
starting from an arbitrary affine subspace $\cal{D}\subseteq\pre(P)$,
containing the frame structure $\D$. The minimal choice for $\cal{D}$ is
to contain only $\D$. In this case the calculus is determined only by
the curvature map $\varrho_\D$. Interestingly, even in this case the
compatibility requirement may imply very strong restrictions on the
calculus. It is also worth noticing that generally the compatibility
between the appropriately chosen subspace $\cal{D}$ will ensure the
full compatibility between $\Gamma$ and $\pre(P)$. 

Applying to $\bigl\{\hor_P,F^\wedge,\pre(P),\Omega_M\bigr\}$
a general construction presented in \cite{d3} we
obtain in the intrinsic way a graded-differential *-algebra
$\Omega(P)$ representing the complete differential
calculus on $P$, such that $\hor_P$ is recovered as the
corresponding *-subalgebra of horizontal forms for $\Omega(P)$.
Moreover, to every map $D$ it is possible to associate intrinsically
a connection $\omega_D$ on $P$, such that $D$ is reinterpreted as
the covariant derivative of the connection $\omega_D$. 

Let us sketch this construction of the calculus $\Omega(P)$. At the level of
graded vector spaces we define
$$
\Omega(P)=\hor_P\otimes\Gamma_{\inv}^\wedge,
$$
while the *-algebra structure is specified by the formulas
\begin{align}
(\psi\otimes\vartheta)(\varphi\otimes\eta)&=
(-1)^{\partial\varphi\partial\vartheta}\sum_k\psi\varphi_k\otimes
(\vartheta\circ c_k)\eta\\
(\varphi\otimes\vartheta)^*&=\sum_k\varphi_k^*\otimes(\vartheta^*\circ
c_k^*),
\end{align}
where $\Sum_k\varphi_k\otimes c_k=F^\wedge(\varphi)$. We have assumed that
the complete calculus on the structure group $G$ is based on the universal
differential \cite{d1} envelope $\Gamma^\wedge$ of $\Gamma$.
Another intrinsic choice for the higher-order calculus is the
braided exterior \cite{w2} algebra $\Gamma^\vee$, associated to the natural
braiding $\sigma\colon\Gamma_{\inv}^{\otimes 2}
\rightarrow\Gamma_{\inv}^{\otimes 2}$.

The action of the differential $d\colon\Omega(P)\rightarrow\Omega(P)$ is
given by
\begin{equation}
d(\varphi)=\D(\varphi)+(-1)^{\partial\varphi}
\sum_k\varphi_k\otimes\pi(c_k)\qquad
d(\vartheta)=d^\wedge(\vartheta)+\varrho_\D(\vartheta),
\end{equation}
and extended to the whole $\Omega(P)$ by requiring the graded Leibniz
rule. In the above formulas $\varphi\in\hor_P$ and $\vartheta\in
\Gamma_{\inv}$, while $d^\wedge\colon\Gamma_{\inv}^\wedge
\rightarrow\Gamma_{\inv}^\wedge$ is the corresponding
differential. We have $d^2=0$ and $d*=*d$. Furthermore, the formulas
\begin{equation}
\widehat{F}(\varphi)=F^\wedge(\varphi)\qquad\widehat{F}(\vartheta)=
\adj(\vartheta)+1\otimes\vartheta
\end{equation}
consistently and uniquely determine a homomorphism
$\widehat{F}\colon\Omega(P)\rightarrow\Omega(P)\grten\Gamma^\wedge$. It
turns out that this map is hermitian, and commutes with the corresponding
differentials. Therefore we can say \cite{d2} that $\Omega(P)$ represents a
differential calculus on a quantum principal bundle $P$. The algebra
$\hor_P$ is reconstructed as the horizontal part of $\Omega(P)$, in
other words
$$ \hor_P=\widehat{F}^{-1}\bigl(\Omega(P)\otimes\cal{A}\bigr). $$

Let us consider an arbitrary 
$D\in\pre(D)$ and let us define $E=D-\D$. The formula
\begin{equation}
\omega_D(\vartheta)=1\otimes\vartheta+\chi_E(\vartheta)
\end{equation}
determines a special connection
$\omega_D\colon\Gamma_{\inv}\rightarrow\Omega(P)$ on $P$. This connection
is regular and multiplicative \cite{d2}. Moreover, the
corresponding operators of covariant derivative and curvature coincide
with $D$ and $\varrho_D$ respectively.

\section{Torsion Operators}

In this section we shall introduce general
torsion operators, and analyze their algebraic properties. We shall also
discuss the question of the uniqueness of the Levi-Civita connection.

Let us start from a framed principal bundle $P$, and assume that the algebra
$\hor_P$ is included in the complete calculus $\Omega(P)$, as described in
the previous section.

Let $\omega\colon\Gamma_{\inv}\rightarrow\Omega(P)$ be an arbitrary
connection \cite{d2} on $P$. 
By definition, this means that $\omega$ is a first-order hermitian linear
map satisfying
$$ \widehat{F}[\omega(\vartheta)]=\sum_k\omega(\vartheta_k)\otimes c_k+
1\otimes\vartheta, $$
where $\Sum_k\vartheta_k\otimes c_k=\adj(\vartheta)$.

Furthermore, the covariant derivative $D_\omega\colon
\hor_P\rightarrow\hor_P$ of the connection $\omega$ is given by equality
$$
D_\omega(\varphi)=d\varphi-\sum_k\varphi_k\omega\pi(c_k),
$$
where $\Sum_k\varphi_k\otimes c_k=F^\wedge(\varphi)$. As we have already
mentioned, particularly interesting are {\it regular connections}. They are
characterized by the fact that $D_\omega$ satisfies the graded Leibniz rule. 

If the calculus $\Gamma$ is maximally compatible with the structure of
$\hor_P$ (as the calculus described in the previous section), then
we have the following correspondence,
\begin{equation*}
\Bigl\{\mbox{Regular connections on $P$}\Bigr\}\leftrightarrow
\left\{\begin{gathered}\mbox{Hermitian first-order covariant}\\
\mbox{antiderivations on $\hor_P$}\\
\mbox{extending $\dM\colon\Omega_M\rightarrow\Omega_M$}
\end{gathered}\right\}=\pre(P)
\end{equation*}
induced by the covariant derivative map. 

The curvature of an arbitrary connection $\omega$ is defined by
$$
R_\omega=d\omega-\langle\omega,\omega\rangle.
$$
This is the analog of the classical first structure equation.
Here $\langle\,\rangle$ are the brackets associated \cite{d2}
to an arbitrary embedded differential map $\delta\colon\Gamma_{\inv}
\rightarrow\Gamma_{\inv}\otimes\Gamma_{\inv}$. Explicitly, $\delta$ is
a hermitian map intertwining the corresponding adjoint actions
of $G$ and satisfying
$$ \delta(\vartheta)=\sum_k\vartheta_k^1\otimes\vartheta_k^2,\qquad
d^\wedge(\vartheta)=\sum_k\vartheta_k^1\vartheta_k^2. $$

\begin{defn} A linear map $T_\omega\colon \V\rightarrow\hor_P$ defined by
\begin{equation}
T_\omega^i=T_\omega(\theta_i)=D_\omega(1\otimes\theta_i)
\end{equation}
is called {\it the torsion} of $\omega$. 
\end{defn}

Let us also define components of the curvature
$R_\omega^{ij}=R_\omega\pi(u_{ij})$. In the following proposition
we have collected the most important properties of the torsion operators.

\begin{pro} \bla{i} The map $T_\omega$ is hermitian and the diagram
\begin{equation}\label{cov-T}
\begin{CD} \V @>{\mbox{$T_\omega$}}>> \hor_P\\
@V{\mbox{$\v $}}VV  @VV{\mbox{$F^\wedge$}}V\\
\V\otimes\cal{A} @>>{\mbox{$T_\omega\otimes\id$}}> \hor_P\otimes\cal{A}
\end{CD}
\end{equation}
is commutative.

\smallskip
\bla{ii} We have
\begin{equation}\label{DT}
D_\omega T_\omega^i=-\sum_j \theta_jR_\omega^{ji}-\sum_j
T_\omega^j\rho_\omega(u_{ji}),
\end{equation}
where $\rho_\omega\colon\cal{A}\rightarrow\hor_P^2$ is a linear map given by
\begin{equation}
\rho_\omega(a)=\langle\omega,\omega\rangle\pi(a)+\omega\pi(a^{(1)})\omega
\pi(a^{(2)}). 
\end{equation}
The map $\rho_\omega$ is related to the lack of
the multiplicativity of $\omega$, and in particular for regular connections
it vanishes identically. 
\end{pro}

\begin{pf} Equality \eqref{DT} is a direct consequence of a general
expression \cite{d2} for the square of the covariant derivative. The diagram
\eqref{cov-T} and the hermicity of $T_\omega$ follow
from the covariance and hermicity of $D_\omega$. 
\end{pf}

It is worth noticing that \eqref{DT} generalizes the classical second
structure equation.
According to our definition of frame structures, the torsion of
the connection $\omega=\omega_\D$ vanishes. 

In classical geometry, the torsion tensor
(together with the metricity condition) uniquely characterizes every
connection. In particular, the Levi-Civita connection is uniquely
characterized by the vanishing torsion condition. 

In our general framework such a characterization does not longer hold. Let us
analyze regular connections having the same torsion. Such connections are
grouped into real affine subspaces of $\pre(P)$,
which are the orbits under the natural
action of the real vector subspace $\cal{X}\subseteq\prv(P)$ consisting of
linear maps $\chi\colon\cal{A}\rightarrow\hor_P$ satisfying
the additional condition
\begin{equation}\label{sym-12}
\sum_j \theta_j\chi(u_{ji})=0. 
\end{equation}
The above condition is a variant of braided-symmetricity of two first indexes
in the coefficients $\chi_{kji}$ given by
$\Sum_k \theta_k\chi_{kji}=\chi(u_{ji})$.
If the braiding $\tau$ is sufficiently `regular' then this symmetricity
condition together with the `antisymmetricity' condition \eqref{k-antisym}
will imply that $\chi(u_{ij})=0$. According to Lemma~\ref{lem:reduct},
this means that $\chi$ vanishes identically.

\begin{defn} We say that a frame structure on $P$ is {\it regular} iff
$\cal{X}=\{0\}$. In this case regular connections are completely determined
by their torsions, and in particular $\D$
is the unique regular connection with the vanishing torsion.
\end{defn}

If we pass to a more general framework allowing non-compact structure groups
and resigning from the unitarity assumption for $\v \colon\V\rightarrow
\V\otimes\cal{A}$ then, generally, the space $\cal{X}$ will be `large'.
This is in a complete agreement with classical geometry. 

\section{Explicit Coordinate Expressions}

Let us first observe that the following natural decomposition holds
\begin{equation}
\hor_P\leftrightarrow\cal{B}\otimes_{\cal{V}}\Omega_M\leftrightarrow
\Omega_M\otimes_{\cal{V}}\cal{B},
\end{equation}
induced by the product map. It follows that the action of
every covariant derivative map $D\in\pre(P)$ is completely determined
by the restriction on $\cal{B}$. In particular, 
the action on the standard coordinate forms is given by
\begin{equation}
D(\theta_i)=\sum_{\alpha}D(b_{\alpha i})\dM(f_\alpha). 
\end{equation}

Therefore we can introduce the `coordinate' description of the covariant
derivatives, by the formula
\begin{equation}
D(b)=\sum_i X_i(b)\otimes \theta_i,
\end{equation}
where the maps
$X_i\colon\cal{B}\rightarrow\cal{B}$ are counterparts of horizontal
coordinate vectors fields. They completely determine
the map $D$. A particular role is given to the
vector fields $\partial_i$ corresponding to $\D$. We have
\begin{equation}
\D(b\otimes\vartheta)=\sum_i\partial_i(b)\otimes\theta_i\vartheta. 
\end{equation}

We are going to analyze how elementary properties of $D$ are reflected on
the maps $X_i$. At
first, $D$ acts on $\cal{B}$ as a derivation. This implies
\begin{equation}\label{X-der}
X_i(qb)=qX_i(b)+\sum_{jk}\mu_{ji}(c_k)X_j(q)b_k,
\end{equation}
where $\Sum_kb_k\otimes c_k=F(b)$ and maps
$\mu_{ji}\colon\cal{A}\rightarrow\Bbb{C}$ are given by
\begin{equation}
\sum_i\mu_{ji}(a)\theta_i=\theta_j\circ a.
\end{equation}

The fact that every covariant derivative $D$ extends
the differential $\dM\colon\Omega_M\rightarrow\Omega_M$ imples
\begin{equation}\label{D-dM}
(X_i-\partial_i)(\cal{V})=\{0\}. 
\end{equation}

Further, the covariance of $D$ is equivalent to equalities
\begin{equation}
FX_i(b)=\sum_{jk} X_j(b_k)\otimes c_k\k^{-1}(u_{ij}). 
\end{equation}

Conversely, let us assume that a system of maps
$X_i\colon\cal{B}\rightarrow\cal{B}$ is given,
satisfying all the above equations. Then the formula
\begin{equation}
D(bw)=\sum_i\bigl(X_i(b)\otimes\theta_i\bigr)w+b\dM(w)
\end{equation}
consistently defines a first-order
linear map $D\colon\hor_P\rightarrow\hor_P$.
This map intertwines $F^\wedge$, extends $\dM$ and satisfies the
`partial' Leibniz rules
\begin{align*}
D(\varphi w)&=D(\varphi)w+(-1)^{\partial\varphi}\varphi \dM(w)\\
D(b\varphi)&=D(b)\varphi+bD(\varphi),
\end{align*}
for each $\varphi\in\hor_P$.

\begin{lem} The following conditions are equivalent:

\bla{i} The map $D$ satisfies the graded Leibniz rule.

\bla{ii} We have
$$ \dM(f)D(b)+D\bigl[\dM(f)b\bigr]=0 $$
for each $f\in\cal{V}$ and $b\in\cal{B}$.

\smallskip
\bla{iii} It is possible to introduce the curvature map $\varrho_D\colon
\cal{A}\rightarrow\hor_P$, naturally associated to $D$. \qed
\end{lem}

Let us now assume that the conditions figuring in the above proposition
are satisfied for $D$. Then the conjugate map $D^*$ is also a covariant
first-order antiderivation extending $\dM$, and we have the unique
decomposition
\begin{equation}
D=D_1+iD_2, \qquad\quad D_1=\frac{D+D^*}{2}\quad D_2=\frac{D-D^*}{2i}. 
\end{equation}
with $D_1,D_2\in\pre(P)$. 
\section{Concluding Examples $\&$ Remarks}
\subsection{Quantum Homogeneous Spaces}

A large class of examples of framed quantum principal bundles
is associated to the appropriate quantum homogeneous spaces.

Let us consider a quantum group $H$, represented by a Hopf *-algebra
$\AH$. Let us assume that $G$ is realized as a subgroup of $H$. This means
that a projection *-homomorphism $\pHG\colon\AH\rightarrow\cal{A}$ is
given such that 
$$
\phi \pHG=(\pHG\otimes \pHG)\phiH\qquad
\e \pHG=\eH \qquad \k\pHG=\pHG\kH.
$$
All primed entities appearing in this subsection refer to the group $H$.

Let $\AHi\subseteq \AH$ be the fixed-point subalgebra relative to
the natural right action $\phiHG\colon\AH\rightarrow\AH\otimes\cal{A}$
of $G$. Furthermore, let us denote by $\adHG\colon\AH\rightarrow
\AH\otimes\cal{A}$ the restricted adjoint action, given by
$\adHG(q)=q^{(2)}\otimes\pHG\bigl[\kH(q^{(1)})q^{(3)}\bigr]$. We have
\begin{align*}
\phiH(\AHi)&\subseteq\AH\otimes\AHi\\
\adHG(\AHi)&\subseteq\AHi\otimes\cal{A}. 
\end{align*}

Let us consider a first-order *-calculus $\Phi$ over $H$, which is
left-covariant and right $G$-covariant. The last condition means that
we can consistently introduce the right action map
$\adjr\colon\Phi_{\inv}\rightarrow\Phi_{\inv}\otimes\cal{A}$, which is
given by the formula
\begin{equation}
\adjr\piH=(\piH\otimes\id)\adHG. 
\end{equation}
It is worth noticing that if $\Phi$ is bicovariant then it is automatically
right $G$-covariant and
$$\adjr=(\id\otimes\pHG)\adjH,\qquad\quad\adjH\colon\Phi_{\inv}\rightarrow
\Phi_{\inv}\otimes\AH.$$

The calculus $\Phi$ is projected via the map $\pHG$, onto the
bicovariant *-calculus $\Gamma$ over the group $G$. Let us denote by
the same symbol $\pHG\colon\Phi_{\inv}\rightarrow
\Gamma_{\inv}$ the corresponding induced map, so that we have
$\pHG(\Phi_{\inv})=\Gamma_{\inv}$,
and 
\begin{equation} \adj \pHG=(\pHG\otimes \id)\adjr\qquad
\pHG(\vartheta\circ q)=\pHG(\vartheta)\circ \pHG(q).
\end{equation}

It follows that $\L=\ker\bigl(\pHG{\restr}\Phi_{\inv}\bigr)$ is a *-submodule of
$\Phi_{\inv}$. Furthermore, we have
$\piH(\AHi)\subseteq \L$. Let us assume
that $\piH\colon\AHi\rightarrow\L$ is surjective.
The space $\L$ is $\adjr$-invariant and we have
\begin{equation}
\adjr\piH(q)=\piH(q^{(2)})\otimes\k\pHG(q^{(1)}), \quad\forall q\in\AHi. 
\end{equation}

Our next assumption is that the right $\AH$-module structure
on $\Phi_{\inv}$ is projectable to the right $\cal{A}$-module structure on
the same space. In other words, 
\begin{equation}
\Phi_{\inv}{\circ}\bigl\{\ker(\pHG)\bigr\}=\{0\},
\end{equation}
and therefore it is natural to
consider $\Phi_{\inv}$ as a right $\cal{A}$-module.

Furthermore, let us assume that the space $\L$ admits a complement
$\L^\perp$ in $\Phi_{\inv}$, which is invariant under the action
$\adjr$, and which is a submodule of $\Phi_{\inv}$. Without a
lack of generality, we may assume that $\L^\perp$ is also *-invariant.
We see that $\L^\perp$ is naturally isomorphic to $\Gamma_{\inv}$. Let
$p_\L\colon\Phi_{\inv}\rightarrow\L$ be the projection map associated to
the spitting $\Phi_{\inv}=\L\oplus\L^\perp$.

The space $\L$, equipped with the induced maps $\{*, \adjr, \circ\}$
determines a bicovariant *-bimodule $\Psi\leftrightarrow\cal{A}\otimes
\L$ over $G$, with the associated braiding $\tau$. 

Finally, let us assume that the implication
\begin{equation}\label{d}
\pHG(q)=a\in\cal{R}\,\Rightarrow\,
\pi_\L(q^{(2)})\pi_\L[\kH(q^{(1)})q^{(3)}]=0
\end{equation}
holds in the corresponding $\tau$-exterior algebra $\L^\wedge$. 
Here $\pi_\L=p_\L\piH$. The above condition is equivalent to a possibility
of introducing consistently a `partial transposed commutator' map
$\Delta\colon\Gamma_{\inv}\rightarrow\L\otimes\L/\im(I+\tau)$ by the formula
\begin{equation}\label{part-com}
\Delta\pi(a)=\pi_\L(q^{(2)})\pi_\L[\kH(q^{(1)})q^{(3)}],\qquad\pHG(q)=a. 
\end{equation}

Let us consider a *-algebra $\cal{B}$, equipped with a free action
$E\colon\cal{B}\rightarrow\cal{B}\otimes\AH$ of $H$. Let $F\colon
\cal{B}\rightarrow\cal{B}\otimes\cal{A}$ be the *-homomorphism given by
$$ F=(\id\otimes\pHG)E. $$
Let $\cal{V}\subseteq\cal{B}$ be the $F$-fixed point subalgebra, and
$i\colon\cal{V}\rightarrow\cal{B}$ the inclusion map. By definition,
$P=(\cal{B},i,F)$ is a quantum principal $G$-bundle over a quantum space
$M$ described by $\cal{V}$. 

Let $\hor_P$ be the horizontal
algebra, associated to $P$ and $\L$, as described at the beginning of
the paper.  Let $\D\colon\hor_P\rightarrow\hor_P$ be a linear
map defined by
\begin{equation}\label{def-D}
\D(b\otimes\vartheta)=\sum_k b_k\otimes\bigl[\pi_\L (q_k)\bigr]\vartheta,
\end{equation}
where $E(b)=\Sum_k b_k\otimes q_k$. 

\begin{lem} The map $\D$ is a hermitian first-order antiderivation
on $\hor_P$, intertwining the map $F^\wedge$. Moreover, there exists the
curvature map associated to $\D$.
\end{lem}

\begin{pf} From the definition~\eqref{def-D} it follows that $\D$ acts
on $\cal{B}$ as a derivation. A direct computation gives
\begin{multline*}
\D(\vartheta b)=\sum_k\D(b_k)\bigl(\vartheta\circ\pHG(q_k)\bigr)=\sum_k\bigl(
b_k\otimes \pi_\L(q_k^{(1)})\bigr)\bigl(\vartheta\circ
\pHG(q_k^{(2)})\bigr)\\
=(-1)^{\partial\vartheta}\sum_k b_k
\otimes\bigl[\bigl(\vartheta\circ
\pHG(q_k^{(1)})\bigr) \pi_\L(q_k^{(2)})\bigr]=
(-1)^{\partial\vartheta}\vartheta \D(b),
\end{multline*}
and we conclude that $\D$ satisfies the graded Leibniz rule. We also see that
$\D$ commutes with the action $F^\wedge$. Furthermore,
\begin{multline*}
\D\bigl[(b\otimes\vartheta)^*\bigr]=\sum_k
\D\bigl[b_k^*\otimes(\vartheta^*\circ
\pHG(q_k^*))\bigr]=\sum_k b_k^*\otimes \bigl[\pi_\L(q_k^{(1)*})(\vartheta^*\circ
\pHG(q_k^{(2)*}))\bigr]\\
=\sum_k (-1)^{\partial\vartheta}\vartheta^*\bigl[
b_k^*\otimes\pi_\L(q_k^*)\bigr]=[\D(b\otimes\vartheta)]^*,
\end{multline*}
which proves the hermicity of $\D$. 

Let us prove that there exists the curvature map. We compute
$$ \D^2(b)=\sum_k b_k \pi_\L(q_k^{(1)})\pi_\L(q_k^{(2)}), $$
and hence if the curvature map exists, it is necessarily given by the
formula
$$ \varrho_\D\pHG(q)=-\pi_\L(q^{(1)})\pi_\L(q^{(2)}).$$
We have to check the consistency of the above formula. However, this is a
direct consequence of the self-consistency of the partial commutator
map $\Delta$, given by \eqref{part-com}. Explicitly, two maps are related
by
\begin{equation}\label{rel:r-d}
\varrho_D(\vartheta)=-\frac{1}{2}\Delta(\vartheta),
\end{equation}
where now $\vartheta\in\Gamma_{\inv}$.
Indeed, we have
\begin{equation}
\tau\bigl(\pi_\L(q^{(1)})\otimes\pi_\L(q^{(2)})\bigr)=
\pi_\L(q^{(1)})\otimes\pi_\L(q^{(2)})-
\pi_\L(q^{(2)})\otimes\pi_\L\bigl(\kH(q^{(1)})q^{(3)}\bigr),
\end{equation}
and \eqref{rel:r-d} follows by projecting down to $\L^\wedge$. 
It is worth observing that \eqref{d} and
\eqref{rel:r-d} prove also the automatical
compatibility between the calculus $\Gamma$ and the curvature $\varrho_\D$. 
\end{pf}

Let us fix a basis $\{\theta_1,\dots,\theta_n\}$ in the space $\L$. 
Let us consider the elements $c_1,\dots, c_n\in\AHi$ satisfying
$$
\adHG(c_i)=
\sum_j c_j\otimes u_{ji}\qquad \piH(c_i)=\theta_i,
$$
where $u_{ji}$ are the matrix elements of
$\v\colon\L\rightarrow\L\otimes\cal{A}$.
Let us define the elements $b_{\alpha i}\in\cal{B}$ and $f_\alpha\in\cal{V}$
by equalities
$$
\sum_\alpha b_{\alpha i}E(f_\alpha)=1\otimes c_i\qquad F(b_{\alpha i})=
\sum_j b_{\alpha j}\otimes u_{ji}. 
$$
A direct calculation shows that \eqref{D-coor} holds, and hence

\begin{pro}
The map $\D$ determines a frame structure on a quantum principal
bundle $P$. \qed
\end{pro}

As a concrete class of examples of the presented construction, let us mention
standard quantum homogeneous spaces (equipped with the appropriate
differential structure $\Phi$). In this case $\cal{B}=\AH$ and
$E=\phiH$. Another generic class of examples is based on quantum classifying
spaces~\cite{d-qclsp}, where we take $\cal{B}$ to be the universal *-algebra
generated by elements $\psi_{ki}$, where $k\in\{1,\dots,d\}$ and
$i\in\{1,\dots n\}$, modulo the relations $\Sum_k\psi_{ki}^*\psi_{kj}
=\delta_{ij}1$. The base space $M$ is one of the classifying spaces for $G$,
and $P$ is the analog of the universal $G$-bundle. Differential calculus
constructed using the presented methods is particularly suitable~\cite{d-cl}
for constructing examples of quantum characteristic classes.

\subsection{Classical Structure Groups}

A specially interesting class of examples is given by quantum frame bundles
with classical~\cite{d-frm} structure groups $G$. In this context we can
further assume that the right $\cal{A}$-module structure $\circ$ on $\V$
is trivial, and that the frame structure $\D$ is compatible with the
classical differential calculus on $G$. This is equivalent to
\begin{align}
\vartheta\circ a&=\e(a)\vartheta\label{cl-mod}\\
\varrho_\D(ab)&=\e(b)\varrho_\D(a)+\e(a)\varrho_\D(b)\label{cl-rD}
\end{align}
respectively, and in particular
all flip-over operators become the standard transpositions.

If we consider $G=\SO(n)$ and assume that $u$ is the standard
representation of $G$ in $\V=\Bbb{C}^n$,
then if $n=2k+1$ the unique 
solution for the right-module structure is given by \eqref{cl-mod}. On the other hand,
in the case $n=2k$ we have also a possible solution
\begin{equation}
\vartheta\circ u_{ij}=-\delta_{ij}\vartheta.
\end{equation}
For $k\geq 2$ this will be the unique nontrivial
solution for the
right-module structure. We see that $-\tau$ is the standard
transposition on $\V$ and hence $\V^\wedge$ is consisting of the symmetric
tensors over $\V$. Frame structures based on such
$\{\V,\circ\}$ give interesting variants
of the `fermionic' differential calculi on the base space.

The presented formalism effectively incorporates \cite{d-qkm} 
the basic elements of the classical theory of {K\"a}hler manifolds.
In this context the structure group is reduced to $G=\U(k)$. 

However, even if $\circ$ is trivial and $\D$ is arbitrary the minimal
calculus $\Gamma$ will be non-standard. 

In what follows we shall analyze in more details framed complex
line bundles.
This is specified by $G=\SO(2)$, and the standard $2$-dimensional
representation
\begin{equation*}
\v =\begin{pmatrix} \cos\varphi & -\sin\varphi\\
\sin\varphi &\phantom{-}\cos\varphi
\end{pmatrix}
\end{equation*}
where $\varphi$ is the angle function. Let us assume that the *-structure is
trivial. We have to fix the $\circ$-operation in the space $\V=\Bbb{C}^2$.
The solutions are given by
\begin{equation*}
\begin{aligned}
\psi\circ u&=\lambda\psi\\
\psi^*\circ u&=\lambda\psi^*
\end{aligned}\qquad\quad
\begin{aligned}
\psi^*&=\theta_1+i\theta_2\\
\psi&=\theta_1-i\theta_2
\end{aligned}
\end{equation*}
where $0\neq\lambda\in\Re$. In the basis $\{\psi,\psi^*\}$,
the flip-over operator $\tau$ is given by
\begin{equation*}
\tau=\begin{pmatrix} 
\lambda &0&0&0\\
0&0&\lambda&0\\
0&\lambda^{-1}&0&0\\
0&0&0&\lambda
\end{pmatrix}
\end{equation*}
and relations in the exterior algebra $\V^\wedge$ are given by
$$ \lambda\psi\psi^*=-\psi^*\psi \qquad\quad\psi^2=\psi^{*2}=0, $$
where we have assumed that $\lambda\neq-1$.

Let us consider a quantum principal $G$-bundle $P$, equipped with a
frame structure $\D$. It turns out
that the associated calculus on $G$ will be generally infinite-dimensional
(this is caused by the fact that the algebra $\cal{B}$ may be very
complicated, and the components of the curvature can be in principle
arbitrarily complicated elements of this algebra). However, in some special
cases we can perform a further reduction. For example, let us assume that
$\D$ is such that its
curvature takes values from $\V^\wedge\subseteq\hor_P$. In other words,
$$ \varrho_{\D}(a)=c(a)\psi\psi^*$$
with $c\colon\cal{A}\rightarrow\Bbb{C}$. If $\varrho_{\D}$ is non-vanishing
then the minimal calculus $\Gamma$ compatible with the frame structure $\D$
will be $1$-dimensional, and
$$\vartheta\circ u=\lambda^2\vartheta$$
for each $\vartheta\in\Gamma_{\inv}$. This calculus is based on the ideal
$$\cal{R}=\gen\bigl\{1+\lambda^2-u-\lambda^2u^{-1}\bigr\}. $$

For the end of this subsection, let us illustrate the main structural
elements of the presented theory in the context of quantum line bundles.
We shall assume that the base space $M$ is a classical compact smooth
manifold. However, many of the formulas hold in the full generality.

Let $\Lambda$ be a complex line bundle over $M$, and let $\cal{F}$ be the space
of smooth sections of $\Lambda$. Let $\gamma$ be a given automorphism of the
algebra $\cal{V}=S(M)$. We shall assume that $\cal{F}$ is equipped with a
$\cal{V}$ bimodule structure consisting of the standard multiplication on
the left, and the right $\cal{V}$-module structure specified by
\begin{equation*}
\xi f=\gamma(f)\xi. 
\end{equation*}

Now, we can apply the reconstruction procedure as described in \cite{d-tann},
and construct a quantum principal $G$-bundle $P$ over $M$, starting from
$\{\cal{F},\gamma\}$. The algebra $\cal{B}$ is given by
$$ \cal{B}=\sideset{}{^\oplus}\sum_{n\in\Bbb{Z}}\cal{F}_n$$
where
$\cal{F}_k=\cal{F}^{\otimes k}$ and $\cal{F}_{-k}=\bar{\cal{F}}^{\otimes k}$,
for $k\in\Bbb{N}$. The tensor product is taken over $\cal{V}$,
and $\bar{\cal{F}}$ is the conjugate bimodule. The above definition is
completed by saying $\cal{F}_0=\cal{V}$.

The frame structure is determined by a single map $X\colon\cal{B}\rightarrow
\cal{B}$, via the formula
\begin{equation}
\D(b)=X(b)\otimes\psi+\lambda X^*\ll(b)\otimes\psi^*,
\end{equation}
where $\ll\colon\cal{B}\rightarrow\cal{B}$ is an automorphism given by
$(\ll{\restr}\cal{F}_n)=\lambda^n$. So that we have
$$
\psi b=\ll(b)\psi\qquad \ll(b^*)=\ll^{-1}(b)^*
$$
for each $b\in\cal{B}$. 
The map $X$ should satisfy the completeness axiom, and the following
conditions
\begin{gather}
XX^*=X^*X\quad{\restr}\cal{V}\\
X(\cal{F}_n)\subseteq\cal{F}_{n-1}\\
X(bq)=bX(q)+X(b)\ll(q). 
\end{gather}
It follows that
\begin{gather*}
[X,X^*](bq)=[X,X^*](b)\ll(q)+\ll^{-1}(b)[X,X^*](q)\\
X^*(bq)=X^*(b)q+\ll^{-1}(b)X^*(q).
\end{gather*}

The corresponding frame structure $\nabla$ is always regular. The curvature
map is given by
\begin{equation}
\varrho_{\D}(u^n)=\sum_\alpha q_{\alpha}[X,X^*]\ll(b_{\alpha})\otimes
\psi^*\psi,
\end{equation}
where $n\in\Bbb{Z}$ and the elements $q_{\alpha}\in\cal{F}_{-n}$
and $b_{\alpha}\in\cal{F}_n$ satisfy
$$ \sum_\alpha q_{\alpha}b_{\alpha}=1. $$
A direct computation gives
\begin{align}
\varrho_{\D}(u^n)={}&\Bigl\{\sum_{k=0}^{n-1}\lambda^{2k}\gamma^{-k}
\Bigr\}\varrho_{\D}(u)\\
\varrho_{\D}(u^{-n})=-\lambda^{-2}\gamma&\Bigl\{
\sum_{k=0}^{n-1}\lambda^{-2k}\gamma^{k}
\Bigr\}\varrho_{\D}(u),
\end{align}
for each $n\in\Bbb{N}$. In the above formulas $\gamma$ is trivially extended
to $\Omega_M^2=\cal{V}\otimes\{\psi^*\psi\}$.

Therefore, the minimal calculus $\Gamma$ compatible with a single operator
$\D$ will be generally higher-dimensional, including the universal calculus
as a possible minimal solution.

It is worth noticing that the induced calculus over $M$ is non-classical.
Indeed, we have the following natural decomposition
$$ \Omega_M^1=\Omega_M^+\oplus\Omega_M^-, $$
where
$\Omega_M^+=\cal{F}_{-1}\otimes\{\psi\}$ and $\Omega_M^-=
\cal{F}_{1}\otimes\{\psi^*\}$, 
are the spaces playing the role of differential forms of holomorphic and
antiholomorphic types. Therefore 
in $\Omega_M$ we have non-trivial commutation relations.


\begin{thebibliography}{10}
\bibitem{c} Connes A: {\it Noncommutative Geometry}, Academic Press (1994)
\bibitem{d1}
\mbox{\Dj ur\dj evi\'c M: {\it Geometry of Quantum Principal Bundles I},}
Commun Math Phys {\bf 175} (3) 457--521 (1996)
\bibitem{d2}
\mbox{\Dj ur\dj evi\'c M: {\it Geometry of Quantum Principal Bundles II},}
Preprint QmmP 4/93 Belgrade University, Serbia;  
{\it Extended Version:} Preprint, Instituto de Matematicas, UNAM,  
M\'exico (1994)
\bibitem{d3}\Dj ur\dj evi\'c M: {\it Differential Structures on Quantum
Principal Bundles}, Preprint, Instituto de Matematicas, UNAM, M\'exico (1994)
\bibitem{d-cl}\Dj ur\dj evi\'c M: {\it Characteristic
Classes of Quantum Principal Bundles}, Preprint,
Instituto de Matematicas, UNAM, M\'exico (1995)
\bibitem{d-frm}\Dj ur\dj evi\'c M: {\it Classical Spinor Structures
on Quantum Spaces}, Clifford Algebras and Spinor Structures,
Special Volume, 365--377, Kluwer (1995)
\bibitem{d-qkm}\Dj ur\dj evi\'c M:{\it Quantum {K\"a}hler Manifolds},
In preparation, (1996)
\bibitem{d-qclsp}\Dj ur\dj evi\'c M: {\it Quantum Classifying Spaces $\&$
Universal Quantum Characteristic Classes}, Lectures, Stefan Banach
International Mathematical Center, November `95.
\bibitem{d-tann}\Dj ur\dj evi\'c M: {\it Quantum Principal Bundles $\&$
Tannaka-Krein Duality Theory}, Rep Math Phys (to appear, 1996)
\bibitem{ph} Hajac P: {\it Strong Connections and $U_q(2)$ Yang-Mills
Theory on Quantum Principal Bundles}, Preprint, University of California, 
Berkeley (1994)
\bibitem{kn} Kobayashi S, Nomizu K: {\it Foundations of Differential Geometry},
Interscience Publishers, New York, London (1963)
\bibitem{w1} Woronowicz~S~L:{\it Compact Matrix Pseudogroups}, Commun
Math Phys {\bf 111} 613--665 (1987)
\bibitem{w2} Woronowicz~S~L:{\it Differential Calculus on Compact Matrix
Pseudogroups $($Quantum Groups$)$} Commun Math Phys {\bf 122} 125--170 (1989)
\end{thebibliography}
\end{document}